\documentclass[12pt]{article}
\usepackage{graphicx}
\usepackage{tikz}
\usepackage{amssymb}
\title{On a class of meromorphic functions with three singular values}
\author{Alexandre Eremenko}
\def\R{\mathbf{R}}
\def\bR{\mathbf{\overline{R}}}
\def\H{\mathbf{H}}
\def\bH{\mathbf{\overline{H}}}
\def\C{\mathbf{C}}
\def\Z{\mathbf{Z}}
\def\U{\mathbf{U}}
\def\bU{\mathbf{\overline{U}}}
\def\bC{\mathbf{\overline{C}}}

\def\sing{\mathrm{sing}}
\def\diam{\mathrm{diam}\,}

\def\Ima{\mathrm{Im}\, }
\begin{document}
\maketitle
\begin{abstract}
We give a parametrization of the class of real
meromorphic functions in the plane for which the preimage of $\{0,1,\infty\}$
is contained in the real line, and all singular
values belong to $\{0,1,\infty\}$. As corollaries, we obtain a new
parabolic type criterion and prove that
all such functions are limits of rational functions of the same class.

2010 MSC class.: 30D30, Keywords: meromorphic function,
Farey
tessellation, modular function, conformal type.
\end{abstract}

Class $A$ consists by definition of those non-constant
meromorphic functions $f$ in the
complex plane $\C$
which 
have the following properties:
\begin{equation}\label{1}
f(\R)\subset\bR:=\R\cup\{\infty\},
\end{equation}
\begin{equation}\label{2}
f^{-1}(\{0,1,\infty\})\subset\R,
\end{equation}
\begin{equation}\label{3}
\sing(f)\subset \{0,1,\infty\},
\end{equation}
where $\sing(f)$ is the set of all singular values,
that is, critical and asymptotic values\footnote{
When $f:X\to Y$
is a holomorphic map between Riemann surfaces $X$ and $Y$, where $X$ is open,
a curve $\gamma(t),\; 0\leq t<1$, in $X$ is called an asymptotic curve,
if it tends to $\infty$ as $t\to 1$ (in the one-point compactification of $X$)
but $f(\gamma(t))$ tends to some $a\in Y$. This $a$ is called an
asymptotic value. Rational functions are considered on the Riemann sphere $\bC$,
so they do not have asymptotic values.}
of $f$. Transcendental meromorphic functions with
three singular values have been studied recently in \cite{L}, \cite{E},
\cite{BR}. Rational functions with singular values $\{0,1,\infty\}$
are called Belyi functions, and they are relevant for algebraic number theory
\cite{Schneps}.
The Speiser class, consisting of meromorphic functions
with finitely many singular values, plays
an important role in holomorphic dynamics \cite{B} and in the analytic theory
of linear differential equations \cite{BE2}.
\vspace{.1in}

\noindent
{\bf Remark.} It was proved in \cite{BE} that condition (\ref{2}) implies
(\ref{1}), except for an explicitly described family of elementary functions.
\vspace{.1in}

Some examples of functions of class $A$ are 
\begin{equation}
\label{0}
f_1(z)=(\cos z+1)/2,\quad f_2(z)=(\cos\sqrt{z}+1)/2.
\end{equation}
$f_2$ is a special case of a ``periodic end''; see, for example,
\cite[Chap. VI]{GO}.

The definition of class $A$ was
inspired by the papers \cite{HZ} and \cite{T}
which
answer a question from \cite{EKS} by constructing an example of a
meromorphic function satisfying (\ref{1}) and (\ref{2}) whose restriction
to the upper half-plane $\H$ is not of bounded type, that is, not
a quotient of bounded holomorphic functions
in $\H$. Functions constructed in \cite{HZ} and \cite{T} have property
(\ref{3}) as well.

In this paper we give a simple parametric description of class $A$.
The class
is evidently invariant under the action of the real affine group
$z\mapsto az+b,\; a\in\R^*,\; b\in\R$ by pre-composition,
and the group of linear-fractional transformations permuting $0,1,\infty$
acting by post-composition. The subclass of rational functions in $A$
is also invariant with respect to the action of $PGL(2,\R)$
acting by pre-composition.
\vspace{.1in}

\noindent
{\bf 1. The 
Farey tessellation
and the Bethe lattice.}
\vspace{.1in}

In the unit disk $\U$, consider
the inscribed circular triangle with all zero angles and vertices $\{1,i,-1\}$.
By reflecting in the sides and in the sides of reflected triangles we obtain
a tessellation of $\U$ by triangles with zero angles, which is called
the {\em Farey tessellation}
\footnote{G. H. Hardy wrote in Mathematician's Apology that
John Farey ``is
immortal because he failed to understand a theorem which Haros
had proved perfectly fourteen years before''.}
. In what follows, by ``triangle'' we mean one of
the triangles of this tessellation. The sides of triangles are open arcs
of circles, they are disjoint and their endpoints are rational
points on the unit circle. At each rational point infinitely many sides meet.
The word ``arc'' will refer to a side of a triangle of the tessellation.
It is easy to see that Euclidean diameters of triangles tend to $0$,
in the sense that
for every $\epsilon>0$ there are only finitely many triangles of diameter $>\epsilon$. The same applies to arcs.

To this tessellation we associate a tree $T_3$,
which is called the {\em Bethe lattice}
in the physics literature. The vertices of the tree\footnote{The word ``vertex'' can mean a vertex of a polygon or a vertex of a graph,
we hope the meaning will be clear from context.} are triangles of the
tessellation, and two vertices are connected by an edge iff the
corresponding triangles have a common side. All vertices of
$T_3$ have degree $3$. We will consider subtrees of $T_3$
as topological objects: two such trees are equivalent if they
are obtained from
each other by a homeomorphism of $\U$. The word ``edge''
refers to an edge of $T_3$, a topological object,
while the words ``arc'' and ``side'' refer to
a side of the tessellation, a geometric object,
though there is a canonical bijection between
edges and sides. 

An infinite {\em path} is a subgraph of $T_3$ of the form $v_0-v_1-v_2,\ldots$
with infinitely
many vertices enumerated by non-negative integers, with one edge between each
$v_j$ and $v_{j+1}$, for each $j\geq 0$. Since we deal only with infinite
paths in this paper we will omit the word ``infinite''.
 Two paths $(v_j)$ and $(v_j^\prime)$
are called {\em equivalent} if their tails coincide,
that is, there exist
natural integers $N$ and $k$ such that $v_j=v^\prime_{j+k}$
for  $j\geq N$. An equivalence class of paths is called an {\em end} of $T_3$.

A path is oriented in the natural order ($v_{j+1}$ follows $v_j$), and each
vertex $v_j$ of a path with $j\geq 1$ is incident to two edges of the path.
Since all vertices have degree $3$ in $T_3$, there is a third edge of $T_3$
attached
to such a vertex $v_j$ which we call an {\em additional edge}. An additional
edge can point to the right or to the left of the oriented path.

For every vertex $v$ of $T_3$ and every end, there is a unique path with
$v_0=v$ representing this end.

To each end, a unique point in $\partial\U$ is associated. Indeed, the arcs
$C_j$ corresponding to the edges of a path are nested, in the sense that
$C_j$ separates in $\U$ all $C_i$ from all $C_k$ for $i<j<k$,  and diameters
of $C_j$ tend to $0$. So there is a unique limit of $C_j$ as $j\to\infty$.

We call this point the {\em limit of the path} or the limit of the end.
For ends with irrational limit,
the correspondence between ends and their limits is
a bijection onto the set of all irrational numbers in $\partial\U$.

For ends with  rational limits, this correspondence is $2$-to-$1$.
More precisely,
a path has a rational limit iff all but finitely many additional edges
of this path point in the same direction.
Now consider a union of two paths with rational
limits, with a common starting vertex $v_0$ and no other common vertices:
$$\ldots-v_{-1}-v_0-v_1-v_2-\ldots$$ 
and introduce the orientation of this doubly infinite path
according to the ordering of integers.
The rational
limits of these two paths are the same iff all additional edges of the doubly
infinite path point in the same direction, right or left. (This means that for
the original paths the additional edges point in the opposite directions with
respect to their native orientations.)

Let $\mu$ be the modular function. It is obtained by mapping one triangle
of the tessellation conformally onto $\H$, sending the vertices of the triangle
to $0,1,\infty$,
and then using the Schwarz reflection
principle to perform an analytic continuation to the whole unit disk.
This function performs a universal covering
$$\mu:\U\to\bC\backslash\{0,1,\infty\}.$$

\noindent
{\bf 2. Representation of $f$ of the class $A$.}
\vspace{.1in}

\noindent
{\bf Lemma 1.} {\em Every function $f\in A$ has a representation
\begin{equation}\label{4}
f(z)=\mu\circ\phi(z),\quad z\in\H,
\end{equation}
where $\phi:\H\to\U$ is univalent.}
\vspace{.1in}

{\em Proof.} Assume that $f$ is transcendental.
It follows from (\ref{3}) that
$$f:\C\backslash f^{-1}(\{0,1,\infty\})\to \bC\backslash\{0,1,\infty\}$$
is an (unramified) covering, see, for example \cite{Surv} (or \cite{M},
if the reader prefers a primary source).
As a covering over the triply punctured sphere,
$\C\backslash f^{-1}(\{0,1,\infty\})$
is hyperbolic. Let
$$\pi: \H\to \C\backslash f^{-1}(\{0,1,\infty\})$$
be a universal covering. Then $f\circ\pi:\H\to \bC\backslash\{0,1,\infty\}$
is a covering, and since
$\H$ is simply connected it is universal. Since $\mu:\U\to \bC\backslash\{0,1,\infty\}$ is another
universal covering, we conclude that
$$f\circ\pi=\mu\circ \alpha,$$
where $\alpha:\H\to \U$ is an isomorphism.
Since $f^{-1}(\{0,1,\infty\})\subset\R$, we conclude
that $\C\backslash f^{-1}(\{0,1,\infty\})$
contains $\H$; let
$\iota:\H\to \C\backslash f^{-1}(\{0,1,\infty\})$ be the inclusion map,
$\iota(z)=z,\;z\in\H$.
Since $\pi$ is a universal covering, we can lift $\iota$ as
$$\iota=\pi\circ h,$$
where $h:\H\to\H$ is a holomorphic local homeomorphism.

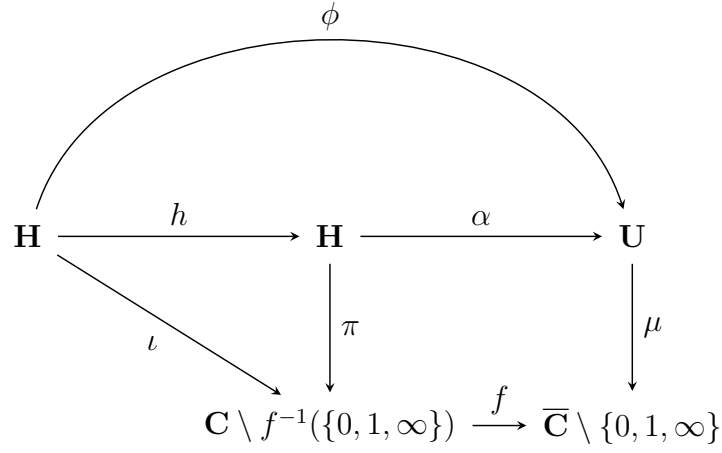
\begin{figure}[htbp]
\centering
\begin{tikzpicture}[>=stealth, line width=0.55pt,
  every node/.style={font=\normalsize},
  map/.style={->,shorten >=2pt,shorten <=2pt}]
  \node (H1) at (0,2.5) {$\H$};
  \node (H2) at (4.0,2.5) {$\H$};
  \node (U)  at (8.0,2.5) {$\U$};
  \node (X)  at (4.0,0) {$\C\setminus f^{-1}(\{0,1,\infty\})$};
  \node (Y)  at (8.0,0) {$\overline{\C}\setminus\{0,1,\infty\}$};
  \draw[map] (H1) -- node[above] {$h$} (H2);
  \draw[map] (H2) -- node[above] {$\alpha$} (U);
  \draw[map] (H1) to[out=72,in=108,looseness=1.05]
     node[above] {$\phi$} (U);
  \draw[map] (H1) -- node[below left] {$\iota$} (X);
  \draw[map] (H2) -- node[right] {$\pi$} (X);
  \draw[map] (U) -- node[right] {$\mu$} (Y);
  \draw[map] (X) -- node[above] {$f$} (Y);
\end{tikzpicture}
\caption{The commutative diagram used in the proof of Lemma~1.}
\label{fig:lemma1-diagram}
\end{figure}

The diagram of the maps we consider is shown in Fig.~\ref{fig:lemma1-diagram}.

Now $h$ is injective since it is a right factor of an injective map $\iota$.
So
we define $\phi=\alpha\circ h$, and obtain (\ref{4}) with an injective
$\phi:\H\to\U$.

The case of rational $f$ is treated similarly;
we leave the details to the reader.
$\Box$
\vspace{.1in}

Let $\Omega:=\phi(\H)\subset \U$ and $\bH=\H\cup\bR=\H\cup\R\cup\{\infty\}$. 
\vspace{.1in}

\noindent
{\bf Lemma 2.} {\em The image $\Omega$ of $\phi$ is a Jordan region in $\C$,
and thus $\phi$ has a homeomorphic extension to
$\bH\to\overline{\Omega}$.
There are three possible types of $\Omega$:
\vspace{.1in}

\noindent
(i) $\Omega$ consists of finitely many triangles and arcs;
the boundary $\partial\Omega$ consists of finitely many arcs and
finitely many points in $\bU$.
In this case $f$ is rational, and $\deg f$ is the number
of triangles in $\Omega$.
\vspace{.1in}

\noindent
(ii) $\Omega$ consists of infinitely many triangles and arcs; the boundary
$\partial\Omega$ consists of an infinite chain of
arcs $C_0,C_1,C_2,\ldots$, and points $y_j\in\partial\U,\; j\geq 0,$ where $C_j$ is tangent
to $C_{j+1}$ and these arcs tend to the endpoint of $C_0$ which is different
from $y_0$.
\vspace{.1in}

\noindent
(iii) $\Omega$ consists of infinitely many triangles and arcs; the boundary
$\partial\Omega$ consists of a doubly infinite chain of arcs
$\ldots C_{-1},C_0,C_1,C_2,\ldots$, points $y_j,\; j\in\Z,$
where $C_j$ is tangent to
$C_{j+1}$, and one more point which is a limit of $C_j$ as $|j|\to\infty$.}
\vspace{.1in}

{\em Proof.} Let $J$ be a component of $\R\backslash \{0,1\}$.
Since both $\mu$ and $f$ are coverings over $\C\backslash\{0,1\}$,
each component of each preimage $f^{-1}(J)$, $\mu^{-1}(J)$ is a simple arc.
For $\mu$ it is an arc $C$ of the tessellation. Let $I$ be a component
of $f^{-1}(J)$ which belongs to the real line. It is clear that
$\phi$ extends to $I$ continuously, and maps it onto some $C$ bijectively.
The endpoints of $I$ must be $a$-points of $f$ with $a\in\{0,1,\infty\}$; 
at these points the limits of $\phi$ must exist, and must be common endpoints
of some arcs $C_j,C_{j+1}$. So we have a chain $\ldots,C_j,C_{j+1},\ldots$,
which can be infinite in one or two directions. This chain together with
common tangent points of its arcs makes a curve $\Gamma$ parametrized by
the real line. Let $y_j\in\partial\U$ be the points where $C_j$ and $C_{j+1}$
are tangent.

To prove that $\gamma=\phi(\R)$ is a simple curve,
we first notice that each arc
$C_j$ separates $\U$ into two regions,
and that $\phi(\Omega)$ must be contained
in one of them, since $\phi$ is univalent.
Next, all arcs $C_j$ are distinct, and thus disjoint. Indeed, if two different
intervals of the real line are mapped by $\phi$ on the same arc $C$, then
the images of a half-neighborhood of these intervals
in $\H$ lie either on the same or on the different sides of this arc.
The first
possibility contradicts univalence of $\phi$ and the second one contradicts the
previous remark that $\Omega$ lies on one side of each arc $C_j$.

So the arcs $C_j$ are disjoint, and the only possible self-intersections
of the
curve $\gamma$ which consists of these arcs and some isolated real points can
happen on the real line.

For such
a self-intersection, it is necessary that at least $3$
arcs $C_j$ share a common end point.
But this is impossible since $\Omega$ lies on one side of each
of the arcs $C_j$.

It is easy to see that the sequence $y_j$ in $\partial U$ is
monotone for $|j|\geq j_0$. So it has a limit. If it is a doubly
infinite sequence, then it has limits when $j\to+\infty$ and when $j\to-\infty$.

Since $\diam C_j\to 0$, the chain must have the same limit(s) as the
sequence $y_j$; if it is two-sided
then limits exist both on the right and on the left. Each such limit
is an accumulation point of images of $a$-points of $f$ under $\phi$.
Since $\phi$ is univalent, the chain cannot have two different
limits along two directions
of the real line. This follows from a theorem of Lindel\"of applied to 
univalent functions or to functions with positive imaginary part. Lindel\"of's
theorem says that such a function cannot have two
distinct asymptotic values with asymptotic curves tending to the same point.
By a small perturbation,
two different limits along $\R$ as $x\to\pm\infty$
imply the existence of
two asymptotic curves in $\H$ tending
to $\infty$ with different asymptotic values, so Lindel\"of's theorem applies.

Now, as  the right and left
limits of our chain of arcs $C_j$ must
coincide, the simple curve $\gamma$ with one point added to it
is the boundary of some Jordan domain
$\Omega_1\subset \Omega$. To see that in fact $\Omega_1=\Omega$,
let $\phi_1$
be a conformal map $\H\to\Omega_1$, normalized by $\phi_1(i)=\phi(i)$ and
$\phi_1(0)=\phi(0)$,
and consider the composition
$\phi^{-1}\circ\phi_1$. This composition has positive
imaginary part in $\H$ and maps
the real line into itself.
So $\Ima \phi^{-1}\circ\phi_1$ is a positive harmonic function in $\H$
with zero boundary values everywhere on $\R$, so it must be $c\Ima z$.
Thus $\phi^{-1}\circ\phi(z)=cz+a$, and by our normalization
it has to be the identity function.
So $\phi=\phi_1$.

This proves that $\phi(\H)$ is a Jordan domain,
so by Caratheodory's theorem
on the boundary correspondence under conformal maps, $\phi$ extends to
a homeomorphism $\bH\to\overline{\Omega}$. The rest of the statements
of the lemma easily follow. $\Box$
\vspace{.1in}

In what follows we mean by $\phi$ the continuous extension to $\bH$.
\vspace{.1in}

\noindent
{\bf Remark.} Let $m$ be the multiplicity of an $a$-point $x$ of $f$,
where $a\in\{0,1,\infty\}$, and $w=\phi(x)$. Then $w=w_j$ for some $j$ as
in Lemma 2, and a
local consideration shows that $m$
can be interpreted as the number of triangles of the tessellation
which share the vertex at $w_j$ and lie between $C_j$ and $C_{j+1}$. 
\vspace{.1in}

\noindent
{\bf 3. A tree corresponding to $f$ and the main result.}
\vspace{.1in}

Consider the part $T=T_f$
of the tree $T_3$ which corresponds to the elements of the tessellation which
intersect $\Omega=\phi(\H)$.

More precisely: the vertices of $T$ correspond
to those triangles of the tessellation which lie in $\Omega$. Edges of $T$
can be of two types: some edges connect two vertices, they correspond to common
sides of pairs of triangles in $\Omega$. Other edges {\em are incident to only
one vertex}; they correspond to the boundary arcs of $\Omega$.
Such edges are {\em semi-open intervals}.
The edges of this last type will be called {\em leaves} of $T$.

So all vertices of $T$ have degree $3$, while an edge can be incident to
one or two vertices.
The only subtree of $T_3$ which has no leaves is $T_3$ itself.

It is convenient to extend each leaf so that it has a definite end in
$\partial\U$; this makes $T$ properly embedded in $\U$,
but this improper end of a leaf is not counted as a vertex!
Once $T$ is properly embedded, we can talk of the {\em faces} of $T$
which are the components of
$\U\backslash T$. A face is called {\em finite} if its boundary contains
finitely many edges of $T$, and {\em infinite} otherwise.
Finite faces correspond to $a$-points of $f$,
$a\in\{0,1,\infty\}$; the number of edges
on the boundary of a finite face is the
multiplicity of the $a$-point plus $1$.
We will
later see that infinite faces correspond to the tracts over asymptotic values
of $f$.

Now consider an arbitrary connected infinite subtree $T$
of $T_3$ of the described
type: all vertices are of degree $3$, and there can be leaves incident to
one vertex, while all other edges are incident to $2$ vertices.
\vspace{.1in}

\noindent
{\bf Definition.}
{\em Such a tree $T$ is called {\em admissible}
if it contains one or two ends,
and in the case of two ends, both have the same limit.}
\vspace{.1in}

We can construct all possible admissible trees as follows.
Begin with an arbitrary path $v_0-v_1-v_2-\ldots$ in $T_3$. Then attach to every
vertex an arbitrary finite subtree of $T_3$ which can be attached
(that is, contains an additional edge of this vertex; $v_0$ has two additional edges, while all the remaining $v_j$ have one).

In this way, we obtain all admissible trees with one end.

For the admissible trees with two ends,
begin with an arbitrary doubly infinite
chain $\ldots-v_{-1}-v_0-v_1-\ldots$ consisting of distinct vertices,
which has the property that all additional
edges are on the same side (with respect to the orientation induced
by the natural order of subscripts).
This is the same as to say that the limits
of the chain as $j\to +\infty$ and as $j\to-\infty$ coincide.
Then attach to each vertex an
arbitrary finite subtree of $T_3$ which can be attached, that is, contains
the additional edge of this vertex. So in this case, all finite trees attached
to the chain are on the same side of the chain.

Two transcendental functions $f,g$ in $A$ are equivalent if
\begin{equation}
\label{e}
f=\ell\circ g\circ m,
\end{equation}
where $\ell$ is a fractional linear transformation
which permutes ${0,1,\infty}$ and $m$ is a real affine transformation.
Two rational functions in $A$ are equivalent if (\ref{e}) holds
with linear fractional $\ell,m$ where $\ell$ permutes $\{0,1,\infty\}$
and $m$ preserves the circle $\overline{\R}$.

Now we are ready to state and prove our main result.
\vspace{.1in}

\noindent
{\bf Theorem 1.} {\em The construction we just described gives
a bijective correspondence
between equivalence classes of transcendental functions of class $A$ and
classes of admissible trees. It also gives a bijection between the
set of equivalence classes of rational functions and classes of
finite subtrees of $T_3$.}
\vspace{.1in}

Thus the class $A$ is parametrized by classes of
finite or admissible subtrees of $T_3$.

{\em Proof}. 1. Let $f$ be a transcendental function of class $A$.
According to Lemma 1, its restriction to $\H$ has a representation
$f=\mu\circ\phi$ which defines $\Omega=\phi(\H)$ and $T=T_f$ as described
above. The resulting tree $T$ must be admissible:
indeed, it is infinite since $f$
is transcendental, so it has at least one path. Moreover,
each path has a limit
and this limit gives an accumulation of $a$-points of $f$.
Different limits would give different accumulation points since
$\phi:\bH\to\overline{\Omega}$ is a homeomorphism by Lemma 2.
So either there is one end or two ends with the same limit. 
\vspace{.1in}

2. In the opposite direction, starting with an admissible tree $T$,
we construct the region $\Omega$. The condition that the tree $T$ is admissible
means that the vertices of triangles belonging to
$\Omega$ have only one accumulation
point; let us call this point $w_\infty$. Now let $\phi:\H\to\Omega$
be the conformal map normalized by $\phi(\infty)=w_\infty$.
Then $\mu\circ\phi$ is real on the real line,
and by reflection extends to a
meromorphic function $f$.
\vspace{.1in}

The rational case is also clear. $\Box$
\vspace{.1in}

Let us explain the interpretation of various features of the tree $T$.
Recall that we extended the leaves to $\partial\U$. 
So the notion of face is
well-defined: it is a component of the complement of $T$. 
Finite faces correspond to $0,1,\infty$ points of $f$.
The number of boundary edges of a
face is equal to the multiplicity of the corresponding $a$-point plus $1$. 
Infinite faces correspond to tracts.
\vspace{.1in}

\noindent
{\bf Corollary.} A function $f\in A$ can have zero, one, or two tracts, but
at most one asymptotic value. Rational functions and 
functions with one end with
an irrational limit have no tracts and no asymptotic values.
A function with one end with a rational limit has one asymptotic value.
A function with two ends has two tracts with the same asymptotic value.
\vspace{.1in}

\noindent
{\bf 4. Line complexes.}
\vspace{.1in}

For general considerations of line complexes
(a. k. a. Speiser graphs or Strecke Complex) and terminology, see
\cite[Chap. VI]{GO}. For a function of class $A$,
choose the extended real line $\bR$ as the base curve.
We have a cell decomposition $C$ of the Riemann sphere with two $2$-cells
(the upper and lower half-planes), three $0$-cells, $0,1,\infty$ and three
$1$-cells $(-\infty,0)$, $(0,1)$ and $(1,+\infty)$.
The preimage of this cell decomposition under $f$ is called the net of $f$,
and the line complex (a. k. a. Speiser graph)
is the preimage of the dual cell decomposition $C'$.
Vertices of the line complex correspond to components of preimages
of half-planes; they are traditionally marked by $\times$ and $\circ$
as they correspond to $\H$ or to the lower half-plane, respectively.
The $1$-skeleton is a bipartite graph embedded in $\C$ (or to $\bC$ in the case
of a rational function) with all vertices
of degree $3$. The faces are marked by $0,1,\infty$ according to their images.

For class $A$, the line complex is symmetric
with respect to the real line, in
the sense that the symmetry maps $\times$-vertices to $\circ$-vertices.
Each bounded face is symmetric. There are none, or one, or two
unbounded faces. If there is only one, it is symmetric; if there are two,
they are interchanged by the symmetry.

The relation between line complexes and trees $T$ is very simple:
Given an admissible tree $T$, consider its reflection $T^*$.
($T^*$ is the image of $T$ under an orientation-reversing involution of $\U$).
Now identify each leaf of $T$ with the leaf of $T^*$ corresponding to it by
the involution. After this identification the two identified leaves become
an ordinary edge and we obtain a line complex. In the opposite direction,
start with a line complex with the stated properties, and cut it into pieces along
the symmetry line. Each piece is an admissible tree.

All resulting line complexes correspond to simply connected
surfaces spread over
the sphere. One consequence of Theorem 1 is that in the case of infinite
three
all such surfaces are of parabolic type.
This does not seem to follow from any usual conformal type
criterion \cite{W}, \cite{Surv}, \cite{Volk}. The symmetry plays a crucial role
here, as was emphasized by Volkovyskii. 
\vspace{.1in}

\noindent
{\bf 5. Approximation of functions $f\in A$ by rational functions in $A$.}
\vspace{.1in}

\noindent
{\bf Theorem 2.} {\em Every function $f$ of class $A$ can be approximated in
the spherical metric
by rational functions of class $A$
uniformly on compact subsets in $\C$.}
\vspace{.1in}

{\em Proof.} For a transcendental function $f\in A$, consider
the tree $T$ corresponding to $f=\mu\circ\phi$,
and let $\Omega=\phi(\H)$ be
the corresponding region. Let
\begin{equation}\label{n1}
w:=\phi(i)\quad\mbox{and}\quad p=\phi(\infty).
\end{equation}
Let $v_0$ be some vertex of $T$.
We exhaust $T$ by a nested sequence of finite subtrees
$T_n$ containing $v_0$, so that $$v_0\in T_1\subset T_2\subset\ldots\to T.$$ 
These finite trees are constructed as follows.
Let $(v_j)$ be the infinite path
containing $v_0$ which is infinite in one or two directions.
(Such a path exists
and is unique for every $v_0$, but a path which is infinite in two
directions
may have coinciding vertices,
namely $v_i=v_{-i}$ for finitely many indices.)
Then $T_n$ is the component containing $v_0$ remaining after cutting the edge
$(v_n,v_{n+1})$, and $(v_{-n},v_{-n-1})$ if the path is infinite in
two directions.

Using $T_n$, construct the regions $\Omega_n$. Evidently
$$\Omega_n\subset\Omega_{n+1}\to\Omega,\quad n\geq 1,$$
and without loss of generality we assume that $w\in\Omega_n,\; n\geq 1$.
Notice that the boundary $\partial\Omega_n$ consists of an arc
$J_n$ of $\partial\Omega$,
and one or two arcs corresponding to the removed edges, and in the
case of two arcs, their common end point $p$. Let
$K_n$ be the closure of the union
of these arcs, so that
\begin{equation}\label{bdrie}
\partial\Omega_n= K_n\cup J_n\quad\mbox{and}\quad
J_n\cup\{ p\}=\partial\Omega_n\cap\partial\Omega.
\end{equation}
Since the diameters of arcs tend to $0$, we have
\begin{equation}\label{diama}
\diam K_n\to 0\quad\mbox{and}\quad\diam(\partial\Omega\backslash J_n)\to0.
\end{equation}
Let $\phi_n:\H\to\Omega_n$ be conformal maps normalized by 
\begin{equation}\label{n2}
\phi_n(i)=w\quad
\mbox{and}\quad
\phi_n(\infty)\in K_n.
\end{equation}
Now we define the rational functions
$f_n=\mu\circ\phi_n$.

It is sufficient to show that $\phi_n\to\phi,\; n\to\infty$ uniformly
on compact subsets of $\H\cup\R$. 
The proof is concluded by the following
\vspace{.1in}

{\bf Lemma 3.} {\em Let $\Omega$ be a Jordan region, $w\in\Omega$, $p\in\partial\Omega$, and
for $n\geq 1$, let $\Omega_n\subset\Omega$ be Jordan regions satisfying
$$w\in\Omega_1\subset\Omega_2\subset\ldots,$$
and 
$$\bigcup\overline{\Omega}_n\supseteq\overline{\Omega}\backslash\{ p\}.$$
We assume that (\ref{bdrie}) and (\ref{diama}) hold. 
Let $\phi:\H\to\Omega$ be the conformal map normalized as in (\ref{n1})
and 
let $\phi_n:\H\to\Omega_n$ be conformal maps normalized by (\ref{n2}).

 Then $\phi_n\to\phi$ uniformly on compact subsets of 
$\H\cup\R$.}
\vspace{.1in}

{\rm Proof.} 
 We denote by $\omega(z,\alpha,D)$ the harmonic measure of
$\alpha\subset \partial D$ with respect to a region $D$.
Condition (\ref{diama}) implies that
\begin{equation}\label{hto0}
\omega(w,K_n,\Omega_n)\to 0\quad\mbox{and}\quad\omega(w,J_n,\Omega)\to 1.
\end{equation}
(To see this, notice that $\omega(z,K_n,\Omega_n)$ is majorized
by $\omega(z,\partial D_\epsilon,\U\backslash D_\epsilon)$, where $D_\epsilon$
is the disk of radius $\epsilon=2\,\diam K_n$ centered at $p$.)

Consider the functions $g_n=\phi^{-1}\circ\phi_n:\H\to\H$.
Our normalizations imply that $g_n(i)=i$. 
As conformal maps
between Jordan regions, $g_n$ extend to the homeomorphisms of the closures, and
we have $$g_n(I_n)=I_n^\prime\subset\R,
\quad\mbox{where}\quad I_n=\phi_n^{-1}(J_n),\quad
\mbox{and}\quad I_n^\prime=\phi^{-1}(J_n).$$
In view of (\ref{hto0}), $I_n$ tends to the real line,
in the sense that the left end of this interval tends to $-\infty$ and
the right end tends to $+\infty$. By reflection $g_n$ extend to the conformal
maps 
$$g_n:\H\cup\H^*\cup I_n\to G_n\cup G^*_n,$$
where $^*$ means complex conjugation.
The domain of $g_n$ contains arbitrarily large disks centered at the origin
and the family $\{ g_n\}$ is normal on every such disk, as a family of univalent
functions satisfying $g_n(\pm i)=\pm i$. So any limit function
of this family is univalent and non-constant in $\C$, so it is affine.
But it also has two fixed points $\pm i$, so every limit function is 
the identity function. $\Box$ 
\vspace{.1in}

This completes the proof of Theorem 2. $\Box$
\vspace{.1in}

The author thanks Mikhail Sodin for discussions and help during preparation 
of this article.

\vspace{.1in}

{\em Department of Mathematics, Purdue University,

West Lafayette, Indiana 47907 
\end{document}